%% file: apr.tex
\documentclass[a4paper,11pt]{article}
\usepackage{CAR}

\usepackage{microtype}
\usepackage{xspace}
\usepackage{amsfonts, amssymb, amsmath}
\usepackage{mathtools}
\usepackage{textcomp}
\usepackage{xfrac}
\usepackage[makeroom]{cancel}
\usepackage{bm}
\usepackage[mathscr]{eucal}

\usepackage{graphicx}
\usepackage{url}
\usepackage{xcolor}

\definecolor{darkred}{RGB}{160,0,0}
\definecolor{darkblue}{RGB}{0,0,160}

\usepackage{hyperref}
\usepackage[sort]{cleveref}

\usepackage[cachedir=apr, frozencache=true]{minted}
\setminted{xleftmargin=0.5\parindent}
\setminted[macaulay2]{mathescape=true}

\usepackage{stmaryrd}

\newcommand{\CIperp}{\mathrel{\text{$\perp\mkern-10mu\perp$}}}

\usepackage{listofitems}
\newcommand{\CI}[1]{{%
  \setsepchar{{:}/{|}/{,}}
  \ignoreemptyitems
  \readlist*\mylist{#1}
  \ifthenelse{\listlen\mylist[] = 2}{\mylist[2]}{}%
  [\mylist[1,1,1] \CIperp \ifthenelse{\listlen\mylist[1,1] = 2}{\mylist[1,1,2]}{\mylist[1,1,1]}%
  \ifthenelse{\listlen\mylist[1] = 2}{{} \mid \mylist[1,2]}{}]%
}}

\makeatletter
\newcommand{\kbldelim}{(}
\newcommand{\kbrdelim}{)}

\newcommand{\kbrowstyle}{\scriptstyle}
\newcommand{\kbcolstyle}{\scriptstyle}

\newlength{\kbcolsep}
\newlength{\kbrowsep}

\newif\ifkbalignright

\newlength{\br@kwd}
\newlength{\k@bordht}

\newcommand{\kbordermatrix}[1]{%
\begingroup
	\setbox0=\hbox{$\left\kbldelim\right.$}
	\setlength{\br@kwd}{\wd0}
	\setbox\@arstrutbox\hbox{\vrule
		\@height\arraystretch\ht\strutbox
		\@depth\arraystretch\dp\strutbox
		\@width\z@}
	\setlength{\k@bordht}{\kbrowsep}
	\addtolength{\k@bordht}{\ht\@arstrutbox}
	\addtolength{\k@bordht}{\dp\@arstrutbox}
	\m@th
	\def\@kbrowstyle{\kbrowstyle}
	\setbox0=\vbox{%
		\def\cr{\crcr\noalign{\kern\kbrowsep
			\global\let\cr=\endline
			\global\let\@kbrowstyle=\relax}}
		\let\\\@arraycr
		\lineskip\z@skip
		\baselineskip\z@skip
		\dimen0\kbcolsep \advance\dimen0\br@kwd
		\ialign{\tabskip\dimen0
			\kern\arraycolsep\hfil\@arstrut$\kbcolstyle ##$\hfil\kern\arraycolsep&
			\tabskip\z@skip
			\kern\arraycolsep\hfil$\@kbrowstyle ##$\ifkbalignright\relax\else\hfil\fi\kern\arraycolsep&&
			\kern\arraycolsep\hfil$\@kbrowstyle ##$\ifkbalignright\relax\else\hfil\fi\kern\arraycolsep\crcr
			#1\crcr}%
	}%
	\setbox2=\vbox{\unvcopy0 \global\setbox5=\lastbox}
	\loop
		\setbox2=\hbox{\unhbox5 \unskip \global\setbox3=\lastbox}
		\ifhbox3
			\global\setbox5=\box2
			\global\setbox1=\box3
	\repeat
	\setbox2=\hbox{$\kern\wd1\kern\kbcolsep\kern-\arraycolsep
		\left\kbldelim
		\kern-\wd1\kern-\kbcolsep\kern-\br@kwd
	\vcenter{\kern-\k@bordht\vbox{\unvbox0}}
	\right\kbrdelim$}
	\null\vbox{\kern\k@bordht\box2}
	\endgroup
}
\makeatother

\newcommand{\Th}{\textsuperscript{th}\xspace}

\newcommand{\bolden}[1]{\textbf{#1}}

\newcommand{\defas}{\coloneqq}

\usepackage{bbm}

\renewcommand{\CC}[1]{\mathcal{#1}} 
\newcommand{\BB}[1]{\mathbb{#1}}

\newcommand{\SF}[1]{\mathsf{#1}}
\newcommand{\RM}[1]{\mathrm{#1}}
\newcommand{\TT}[1]{\text{\texttt{#1}}}
\newcommand{\SR}[1]{\mathscr{#1}}

\newcommand{\PD}{\RM{PD}}

\newcommand{\SymMat}{\RM{Sym}}

\newcommand{\Set}[1]{\left\{\,#1\,\right\}}

\newcommand{\T}{\SF{T}}

\newcommand\bbr[1]{%
  \ensuremath{\bm{[}}%
  #1%
  \ensuremath{\bm{]}}%
}

\newcommand\pr[1]{
  \setsepchar{:}
  \ignoreemptyitems
  \readlist*\mylist{#1}%
  \ifthenelse{\mylistlen = 2}{\mylist[2]}{}%
  \bbr{\mylist[1]}%
}

\newcommand\apr[1]{
  \setsepchar{:}
  \ignoreemptyitems
  \readlist*\mylist{#1}%
  \ifthenelse{\mylistlen = 2}{\mylist[2]}{}%
  \bbr{\mylist[1]}%
}

\newcommand\Bracket[1]{
  \setsepchar{:}
  \ignoreemptyitems
  \readlist*\mylist{#1}%
  \ifthenelse{\mylistlen = 2}{\mylist[2]}{}%
  \bbr{\mylist[1]}%
}

\newtheorem{theorem}{Theorem}
\newtheorem{lemma}[theorem]{Lemma}
\newtheorem{proposition}[theorem]{Proposition}

\newtheorem{remark}[theorem]{Remark}

\newtheorem{example}[theorem]{Example}
\newtheorem{definition}[theorem]{Definition}

\newtheorem{conjecture}[theorem]{Conjecture}

\newtheorem{problem}[theorem]{Problem}

\theoremstyle{empty}
\newtheorem{genthm}{Theorem}

\theoremstyle{nonumberplain}
\theoremheaderfont{\itshape}
\theorembodyfont{\normalfont}
\theoremseparator{.}
\theoremsymbol{\ensuremath{\blacksquare}}

\newcommand{\Macaulay}[1]{%
  \ifthenelse{\equal{\detokenize{#1}}{\detokenize{2}}}{}%
  {\PackageError{tboege-preprint}{Always use \protect\Macaulay2, never \protect\Macaulay\space alone}{}}%
  \texttt{Macaulay#1}\xspace}

\begin{document}



\noindent\rule{\textwidth}{1pt}

\vspace{3mm}


\Aufsatz
{Algebra in probabilistic reasoning}
{Algebra in probabilistic reasoning}
{Tobias Boege}
{TB}
{Tobias Boege (Aalto University, Finland)}
{{%
  \begin{tikzpicture}[y=0.80pt, x=0.80pt, yscale=-0.7, xscale=0.7]
  \input{pappus.tikz}
  \end{tikzpicture}%
}}
{post@taboege.de}

\begin{otherlanguage}{english}

\vspace{3mm}
\begin{multicols}{2}
\noindent



\Ueberschrift{Probabilistic reasoning meets~synthetic~geometry}{sec:intro}

\Ueberschriftu{Reasoning about relevance}

Probabilistic reasoning is an approach, based on probability theory,
to tasks of \emph{decision making under uncertainty}. Random variables
are used to model our uncertainty about the factors that define what a
good decision is. Decision making (certain or uncertain) is a vast field
with different paradigms. For example, one may seek to make the globally
least bad decision given what is known certainly about the state of the
random variables. In probably approximately correct learning, on
the other hand, one tries to keep the errors small, most of the time.

Another part of probabilistic reasoning is concerned with \emph{conditional
independence}. This is a ternary relation among jointly distributed random
variables which gives information of the following sort: ``Given that
factor $C$ is observed, factor $A$ does not influence factor $B$''.
This relation is denoted by $\CI{A,B|C}$. We can also choose to condition
on multiple observed variables or none at all and recover the familiar
notion of stochastic independence $\CI{A,B}$.
Conditional independence (in the following abbreviated to CI) is an
extension of stochastic independence which can accomodate a~priori knowledge
that we may have about the outcome of a subsystem $C$ of random variables,
such as when $C$ is \emph{controlled} in a random experiment.
CI~reveals essential combinatorial information that can guide decision
making when only incomplete data about the state of a system is available.

As a first exercise in conditional independence, note what this relation
specializes to when $A = B$. Suppose the outcome $C$ is known and $\CI{A,A|C}$
holds. Then the outcome of $A$ has no bearing on itself after we observe $C$.
This is absurd, \emph{unless} $C$ reveals everything there is to know about~$A$,
i.e., $A$ takes only a single value that depends on the value $C$ takes.
This is known as a \emph{functional dependence} (FD) and it implies the
existence of a function $f$ such that $A = f(C)$ as random variables.

CI and FD provide basic qualitative information about dependencies among
the observations made in, say, a random experiment in the sciences. But
they play a role in other disciplines that deal with the representation
or processing of \emph{information}. For example, a database in relational
algebra may be seen as a (large) sample from an unknown discrete probability
distribution. The designer of a database will usually anticipate CI and
FD relations in the data (e.g., your zip code functionally determines
your city). The purpose of various \emph{normal forms} for relational
databases is to eliminate undesirable dependencies because they increase
the risk of inconsistencies in the data after updates. Instead, the
database must be factored into multiple ``tables'' according to the
normal form. In spirit, what these normal forms demand is similar to
the factorization of a rank-1 matrix into an outer product
\[
  M = a b^\T.
\]
If $M$ is the probability matrix of the joint distribution of discrete
random variables $A$ and $B$, then $a$ and $b$ are uniquely determined
up to a scalar as the row and column sums of $M$ and they correspond to
the marginal distributions of $A$ and $B$. (With more than two random
variables, we factor a rank-1 tensor.) This representation of $M$ as
$a b^\T$ makes maintaining the independence of $A$ and $B$ in the joint
distribution automatic under updates to the marginal distributions $a$
and $b$, and it saves~space!

\Ueberschriftu{Synthetic statistics}

The reasoning task attached to CI is that of \emph{conditional independence
inference}: given a boolean formula $\varphi$ whose variables are
CI~statements and a family of probability distributions, decide if $\varphi$
is true for every distribution in the family. By writing boolean formulas
in conjunctive normal form, one can restrict this investigation to disjunctive
clauses written in implication form, such as:
\[
  \CI{A,C|B} \wedge \CI{A,B} \;\Rightarrow \CI{A,C}.
\]
This formula is one half of the \emph{semigraphoid property} and it holds
for all random vectors \cite[Appendix~A.7]{Studeny}.

\begin{figurehere}
\centering
\scalebox{0.45}{%
\begin{tikzpicture}[y=0.80pt, x=0.80pt, yscale=-1, xscale=1, inner sep=0pt, outer sep=0pt]
\input{pappus.tikz}
\end{tikzpicture}}
\caption{Pappus's theorem in the projective plane.}
\label{fig:pappus}
\vskip 1.5em
\end{figurehere}

Its meaning should be intuitively clear: suppose that knowing $B$
makes $A$ and $C$ independent, but also that $B$ has no influence
on $A$; then $A$ and $C$ must be independent even without knowledge
of~$B$. One should be cautious, however, of leaping to ``intuitively
clear'' conclusions in CI~inference because the laws of probability
theory sometimes seem to defy intuition. The implications
\begin{align}
  \tag{$\cancel{\,1\,}$} \CI{A,B}                 &\Rightarrow \CI{A,B|C} \\
  \tag{$\cancel{\,2\,}$} \CI{A,B|C}               &\Rightarrow \CI{A,B}   \\
  \tag{$\cancel{\,3\,}$} \CI{A,B} \wedge \CI{A,C} &\Rightarrow \CI{A,{(B,C)}}
\end{align}
are all wrong. The reader is invited to construct random experiments
which falsify them. 

The semigraphoid property allows us to deduce with absolute confidence
from some CI~assumptions other CI~consequences, no matter what the
underlying probability distribution is. A valid implication is also
called a \emph{CI~axiom} or \emph{inference rule}. More restrictions on
the distributions under consideration make more valid inference rules
available for reasoning. In this article, we will consider multivariate
normal distributions. This class has many favorable properties:
(1)~relatively few parameters are needed to specify a distribution,
(2)~they include classes of popular graphical models and
(3)~conditional independence has an \emph{algebraic} reformulation.
Whenever possible we wish to use algebra in reasoning and benefit from
the exactness of symbolic~methods.

In this article, I want to explain a different point of view on
probabilistic reasoning, in particular CI~inference. It is motivated
by similarities to synthetic geometry \cite{BokowskiSturmfels} which
describes geometric objects in relations of ``special position''
to each other. \Cref{fig:pappus} illustrates Pappus's theorem in the
projective plane over a field. This is an inference rule in synthetic
geometry stating that: if all points on the solid lines are collinear,
then the points on the dashed line must also be collinear.
Instead of geometric objects, in probabilistic reasoning we describe
random variables and their ``special position'' in relation to each
other. Special position in this case is conditional independence and
we wish to obtain rules of reasoning such as Pappus's theorem in this
setting --- which we may call \emph{synthetic statistics}.

\Ueberschrift{The geometry of CI~inference}{sec:inf}

\Ueberschriftu{Algebraic statistics of Gaussian~random~vectors}

We suppose a finite \emph{ground set} $N$ of size $n$ indexing the entries
of a random vector $X = (X_i : i \in N)$. Instead of referring to random
variables, subsequently we refer to their indices. It is customary not
to distinguish between an element $i \in N$ and a singleton subset
$\Set{i} \subseteq N$: both are usually denoted by~$i$. Moreover, the
symbol for set union $I \cup K$ for subsets of $N$ is usually~omitted.
Hence, an expression such as $iK$ means the subset $\Set{i} \cup K \subseteq N$.
Denote by $\SymMat_N(\BB K)$ the affine space of $N \times N$ symmetric
matrices over a field $\BB K \subseteq \BB R$ and by $\PD_N(\BB K)$ the
semialgebraic subset of positive definite matrices. Recall that this is
a full-dimensional, open convex cone and that its boundary is the
hypersurface of singular positive semidefinite matrices.

A regular multivariate normal (``Gaussian'') distribution is determined by
its mean $\mu \in \BB R^N$ and its covariance matrix $\Sigma \in \PD_N(\BB R)$.
Regularity refers to the covariance matrix $\Sigma$ which in general only
needs to be positive semidefinite. While an algebraic theory of conditional
independence can be developed even in the degenerate case, this is more
involved and we stick to the regular Gaussians~here.

The density of a Gaussian random vector $X$ with respect to the standard
Lebesgue measure on $\BB R^N$ is a proper transcendental function depending
on $\mu$ and~$\Sigma$:
\[
  \frac1{\sqrt{(2\pi)^n \det \Sigma}}
    \exp\left(-\frac12 (x-\mu)^T \Sigma^{-1} (x-\mu) \right).
\]
A surprising but basic fact of algebraic statistics is that the conditional
independence relation among the components of $X$ depends only on certain
polynomial expressions in the covariance matrix.

\begin{definition}
A \emph{principal minor} of $\Sigma$ is a subdeterminant $\pr{K:\Sigma}
\defas \det \Sigma_{K,K}$ for some $K \subseteq N$. An \emph{almost-principal
minor} is a subdeterminant of the form $\apr{ij|K:\Sigma} \defas
\det \Sigma_{iK,jK}$ for $ijK \subseteq N$ and $i,j \not\in K$ distinct,
where $iK$ indexes rows and $jK$ columns.
\end{definition}

\begin{genthm}[Sign convention] \label{con:sign}
It~is advantageous for the general CI~theory to view the set $N$ and hence
the rows and columns of our matrices as unordered. To still get a
well-defined sign for its determinant, it only matters which row and
column labels $r, c \in N$ are paired together in the $k$\Th position
from the top-left corner of the matrix. For~principal minors $\pr{K:\Sigma}$
and almost-principal minors $\apr{ij|K:\Sigma}$ we establish the following
convention: in the principal submatrix with respect to $K \subseteq N$,
pair each $k \in K$ with itself, and in the almost-principal submatrix
additionally pair row~$i$ with column~$j$.
\end{genthm}

\begin{lemma} \label{lemma:algci}
A symmetric matrix $\Sigma \in \SymMat_N(\BB R)$ is positive definite
if and only if $\pr{K:\Sigma} > 0$ for all $K \subseteq N$. If $\Sigma$
is the positive definite covariance matrix of the Gaussian random vector
$X$, then the conditional independence $\CI{X_i,X_j|X_K}$ holds if and
only if $\apr{ij|K:\Sigma} = 0$.
\end{lemma}

The crucial ingredient to prove this is in \cite[Proposition~4.1.9]{Sullivant}.
In particular, the CI~relation does not depend on the mean $\mu$ and we
may identify Gaussian distributions with their positive definite covariance
matrices.

\begin{definition}
A \emph{Gaussian CI~model} is a subset of $\PD_N(\BB R)$ which is
given by vanishing and non-vanishing constraints on almost-principal
minors (referred to as the \emph{independence} and the \emph{dependence
assumptions} of the model, respectively).
\end{definition}

\begin{remark}
It should be noted that in geometry it is (linear) \emph{dependence}
of vectors which corresponds to a Zariski-closed condition, whereas in
statistics it is the (conditional) \emph{independence} which is closed.
\end{remark}

Gaussian conditional independence models are subsets of the $\PD$ cone
which are cut out by \bolden{very special} classes of polynomial
constraints. They are all determinantal and, up to the symmetric group
on $N$ acting on the coordinates of $\SymMat_N(\BB R)$, there is
precisely one principal and one almost-principal minor of each degree.

\Ueberschriftu{Conditional independence inference}

Consider the general implication formula $\varphi$:
\[
  \bigwedge_{p=1}^s \CI{i_p,j_p|K_p} \;\Rightarrow\; \bigvee_{q=1}^t \CI{x_q,y_q|Z_q}
\]
involving CI~statements over a fixed ground set~$N$. This formula is a
valid inference rule for Gaussians if and only if every $\Sigma \in
\PD_N(\BB R)$ which satisfies $\apr{i_pj_p|K_p:\Sigma} = 0$ for \emph{all}
$p \in [s]$ also satisfies $\apr{x_qy_q|Z_q:\Sigma} = 0$ for \emph{at
least one} $q \in [t]$.
We associate to $\varphi$ a CI~model $\CC M(\varphi)$ which is defined
by the independence assumptions $\CI{i_p,j_p|K_p}$, $p \in [s]$, and the
dependence assumptions $\neg\CI{x_q,y_q|Z_q}$, $q \in [t]$. This is the
set of \emph{counterexamples} to $\varphi$; the implication is valid if
and only if $\CC M(\varphi) = \emptyset$. Conversely, every CI~model is
the set of counterexamples to a suitable CI~implication formula.

Hence, the CI~inference problem is equivalent to the problem of checking
if a system of independence and dependence assumptions is consistent,
which in turn reduces to checking the feasibility of a semialgebraic set
which is defined by integer polynomials --- the CI assumptions as well
as positive definiteness.

\begin{example} \label{ex:weak}
The \emph{weak transitivity} property of Gaussians over $N = ijk$
states that
\[
  \CI{i,j} \wedge \CI{i,j|k} \Rightarrow \CI{i,k} \vee \CI{j,k}.
\]
To prove this rule algebraically, we first determine the ideal
generated by the assumptions:
\begin{align*}
  \apr{ij|\emptyset:\Sigma}  &= \sigma_{ij}, \\
  \apr{ij|k:\Sigma} &= \sigma_{ij} \sigma_{kk} - \sigma_{ik} \sigma_{jk}.
\end{align*}
These vanishing conditions give rise to a hyperplane and a quadratic
hypersurface in the space of $3 \times 3$ correlation matrices, pictured
in \Cref{fig:weak}. The equation $\sigma_{ik} \sigma_{jk} = 0$ holds
on the intersection, which is the union of two line segments: one for each
of the possible conclusions $\CI{i,k} \vee \CI{j,k}$ of weak transitivity
reasoning. The line segments intersect in the identity matrix which
satisfies both conclusions.
\end{example}

\begin{figurehere}
\centering
\includegraphics[width=0.8\linewidth]{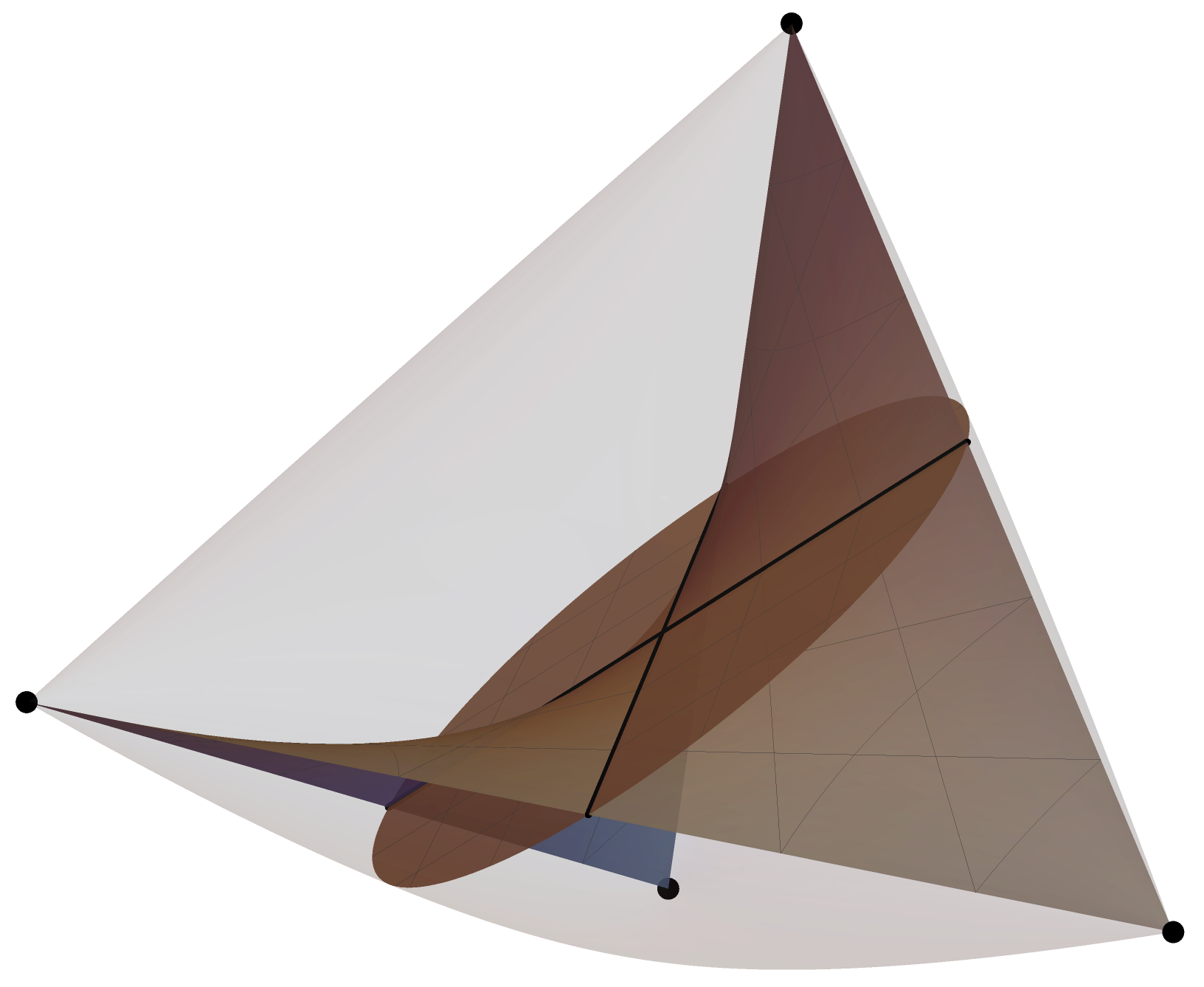}
\caption{Geometric proof of weak transitivity.}
\label{fig:weak}
\vskip 1.5em
\end{figurehere}

It was observed by Matúš \cite{MatusGaussian} that fundamental inference
rules for Gaussians, including a generalization of \Cref{ex:weak} to
arbitrary conditioning sets, follow from a single determinantal identity:

\begin{genthm}[Matúš's identity] \label{lemma:matus}
Let $R$ be a commutative ring with unity and $\Sigma \in \SymMat_N(R)$.
The following identity holds for all $ijkL \subseteq N$:
\begin{gather}
  \label{eq:matus}
  \begin{gathered}
    \pr{kL : \Sigma} \cdot \apr{ij|L : \Sigma} = \qquad\qquad\qquad\qquad\qquad \\
    \qquad \pr{L : \Sigma} \cdot \apr{ij|kL : \Sigma} +
    \apr{ik|L : \Sigma} \cdot \apr{jk|L : \Sigma}.
  \end{gathered}
\end{gather}
\end{genthm}

\begin{remark}
Matúš's formulation of this result in \cite{MatusGaussian} adds a sign
to one of the terms of the identity, depending on the relative ordering
of $i$, $j$ and $k$. This sign is fixed to $+1$ by our Sign convention.
\end{remark}

\begin{remark}
Moreover, Matúš proves this identity over the complex numbers, but the
extension to any commutative ring is standard. It suffices even to prove
the identity only for real positive-definite matrices. Since then the
identity is known to hold on a full-dimensional (and therefore dense)
subset of the irreducible affine space $\SymMat_N(\BB C)$, it holds for
all symmetric matrices over~$\BB C$. But then evaluating the linear
combination of determinants in~\eqref{eq:matus} must result in the zero
polynomial in $\BB Z[\Sigma]$ and thus in every commutative ring with
unity.
\end{remark}

The general weak transitivity property
\[
  \CI{i,j|L} \wedge \CI{i,j|kL} \Rightarrow \CI{i,k|L} \vee \CI{j,k|L}
\]
is a simple consequence of this relation and the algebraic definition
of conditional independence in \Cref{lemma:algci}.

This motivates the investigation of general determinantal identities
which can be formulated in terms of principal and almost-principal
minors only. More precisely, we consider the polynomial ring $\SR R_N$
whose variables are formal brackets $\pr{K}$ and $\apr{ij|K}$ over the
ground set~$N$. (The~empty principal minor bracket $[\emptyset]$ acts
as a homogenization variable.)
Let $\SR J_N$ be the homogeneous ideal given by the kernel of the
evaluation map sending an element of $\SR R_N$ into the coordinate
ring $\BB C[\Sigma]$ of $\SymMat_N(\BB C)$ by evaluating brackets
$\pr{K} \mapsto \pr{K:\Sigma}$ and $\apr{ij|K} \mapsto \apr{ij|K:\Sigma}$.
The generators of the quadratic part of $\SR J_N$ have been found in
\cite{Geometry}.

In some way, the Matús identity, which belongs to this generating set,
acts like the three-term Grassmann--Plücker relations in geometry. While
the latter lead to the definition of \emph{matroid} as a combinatorial
model for geometric special position, the Matúš identity leads to the
\emph{gaussoid axioms}, the most basic inference rules for Gaussian
CI~relations. The gaussoid axioms satisfy two important completeness
results which justify viewing the Matúš identity as fundamental (even
though it alone does not generate the quadratic part of $\SR J_N$).

\begin{proposition}[\cite{LnenickaMatus}]
The gaussoid axioms generate (by logical implication) all true inferences
for three Gaussian random variables.
\end{proposition}

\begin{proposition}[\cite{Gaussant}]
The gaussoid axioms generate (by logical implication) all true inferences,
having at most two assumptions, among any number of Gaussian random variables.
\end{proposition}

The ideal $\SR J_N$ encodes the particular combinatorial flavor of principal
and almost-principal minors of symmetric matrices that distinguishes the
algebraic geometry of Gaussian conditional independence from that of point
configurations and special position which is instead derived from the maximal
minors of rectangular matrices. The following conjecture about $\SR J_N$ is
still open. Its analogue in synthetic geometry is an established theorem
stating that the Grassmann--Plücker relations generate the vanishing ideal
of the Grassmannian.

\begin{conjecture}[\cite{Geometry}]
The ideal $\SR J_N$ is generated by its quadratic part.
\end{conjecture}

\Ueberschrift{Algebraic certificates and~reproducibility}{sec:cert}

There is no finite axiomatization of all inference rules which are valid for
$n$ Gaussian random variables as $n$ grows. The notion of axiomatization can
be made precise by introducing \emph{minors}. In analogy to graph and matroid
theory, minors are the ``natural subconfigurations'' of a collection of
random variables. In graph theory, one obtains minors by deleting and
contracting edges. In probability theory, we take marginal and conditional
distributions. Non-axiomatizability results have been achieved independently
by Šimeček~\cite{SimecekNonfinite} and Sullivant~\cite{Sullivant}.

Nevertheless, one may be interested in finding all valid inference rules for
small numbers of random variables, so that larger systems can be partially
reasoned about, one couple of variables at a time. The case of three Gaussians
is covered by the gaussoid axioms. The four-variate case was solved by
Lněnička and Matúš in~\cite{LnenickaMatus} and the five-variate case is
still wide open.

\columnbreak

This classification task was posed as Challenge~1 in~\cite{Geometry}.
Determining the realizability status of the $254\,826$ candidate gaussoids
computed there is equivalent to finding all valid inference rules on five
Gaussian random variables. A proper solution to such a large-scale
classification task is a FAIR database (cf.~\cite{Mathrepo}) which includes
not only the classification itself but also \emph{machine-checkable proofs}
of its correctness.
The existence of these proofs is guaranteed by theorems in real algebra.
I want to close this exposition by discussing what they are and the
obstacles currently faced when computing them in~practice.

\Ueberschriftu{Algebraic numbers and final polynomials}

Let us return to Pappus's theorem for a moment. Call the three points on
the upper, the lower and the middle line in \Cref{fig:pappus} \TT{a}, \TT{b},
\TT{c}, then \TT{d}, \TT{e}, \TT{f} and then \TT{g}, \TT{h}, \TT{i},
respectively, from left to right. For brevity denote below the determinant
of the $3 \times 3$ matrix those columns are the homogeneous coordinates
of the points labeled $\TT p,\TT q,\TT r$ by the bracket~$\pr{\TT{pqr}}$.
One strikingly mechanic way of proving Pappus's theorem is presented in
the following snippet of \Macaulay2 code:

\begin{minted}[fontsize=\footnotesize]{macaulay2}
-- Homogeneous coordinates for a ... i.
R = QQ[a_1..a_3, b_1..b_3, c_1..c_3,
       d_1..d_3, e_1..e_3, f_1..f_3,
       g_1..g_3, h_1..h_3, i_1..i_3];

-- Bracket is an abbreviation for
-- the 3x3 determinant of (p q r).
br = (p,q,r) -> det matrix(
  apply({p,q,r}, x -> apply({1,2,3},
    i -> value(toString(x)|"_"|i)
  )));

-- The ideal of collinearity assumptions
-- for Pappus's theorem.
papp = ideal(
  br(a,b,c),br(d,e,f),br(a,e,g),br(a,f,h),
  br(b,d,g),br(b,f,i),br(c,d,h),br(c,e,i)
);

-- The conclusion [ghi] and non-degeneracy
-- assumptions.
G = {
  br(g,h,i),br(a,d,i),br(a,b,d),br(a,c,i),
  br(a,d,e),br(a,g,i),br(a,d,h),br(a,f,i),
  br(d,e,i),br(a,d,f),br(d,e,i),br(a,d,f),
  br(d,h,i),br(a,c,d),br(b,d,i),br(a,d,g),
  br(a,b,i),br(d,f,i),br(a,e,i),br(c,d,i)
};

fold((x,y) -> (x*y) % papp, G) --> 0
\end{minted}

This computation producing a zero at the end proves the existence of
polynomials $h_1, \dots, h_8$ with rational coefficients such that
\begin{gather}
  \label{eq:pappus}
  \begin{split}
  &\qquad\qquad[\TT{ghi}] \cdot \prod_{i=2}^{20} g_i = \\
  &h_1 [\TT{abc}] + h_2 [\TT{def}] + h_3 [\TT{aeg}] + h_4 [\TT{afh}] + {} \\
  &h_5 [\TT{bdg}] + h_6 [\TT{bfi}] + h_7 [\TT{cdh}] + h_8 [\TT{cei}],
  \end{split}
\end{gather}
where the polynomials $g_i$ are the elements of the list \mintinline{macaulay2}{G}
in the above code listing. The first polynomial in \mintinline{macaulay2}{G} is the
desired conclusion that \TT{g}, \TT{h} and \TT{i} are collinear. \penalty-1000
All~other polynomials in \mintinline{macaulay2}{G} are non-zero because they
correspond to non-degeneracy conditions for the point and line configuration
in \Cref{fig:pappus}. This implies Pappus's theorem because it shows that
$[\TT{ghi}]$ vanishes, under the non-degeneracy conditions, whenever the
assumptions of Pappus's theorem are satisfied. The crucial list
\mintinline{macaulay2}{G} for this proof is extracted from the fourth
proof of Pappus's theorem in \cite[Section~1.3]{RichterGebert}.

The polynomial $\prod_{i=1}^{20} g_i$ and the linear combinators
$h_1, \dots, h_8$ represent a self-contained proof of Pappus's theorem.
This proof is big and relatively hard to find but verifying it is a matter
of \emph{multiplying out} both sides of \eqref{eq:pappus} and comparing
coefficients. This can be done in any computer algebra system off-the-shelf.
This high standard of verifiability is, fortunately, a theorem which
extends far beyond Pappus:

\begin{genthm}[Positivstellensatz]
The system $\Set{ f_i = 0, g_j \ge 0, h_k \not= 0 }$ defined by finite
collections of polynomials $f_i, g_j, h_k \in \BB Z[x_1, \dots, x_n]$ has
no solution if and only if $0 \in \SR I + \SR P + \SR U^2$, where $\SR I$
is the ideal generated by the $f_i$, $\SR P$ is the cone generated by the
$g_j$ and $\SR U$ is the multiplicative monoid generated by the $h_k$ in
$\BB Z[x_1, \dots, x_n]$.
\end{genthm}

This version of the Positivstellensatz is proved in \cite[Proposition~4.4.1]{RealAlgebra}.
The condition $0 \in \SR I + \SR P + \SR U^2$ implies the existence of an
integer polynomial $f \in \SR I \cap (\SR P + \SR U^2)$. This \emph{final
polynomial} (cf.~\cite[Section~4.2]{BokowskiSturmfels}) must be simultaneously
zero and positive on every point satisfying the polynomial system.
It therefore serves as an obvious proof of the emptiness of the
semialgebraic~set. The coefficients witnessing that $f \in \SR I$
and $f \in \SR P + \SR U^2$ constitute the algebraic certificate.

In the context of probabilistic reasoning, we can now certify when the
model of counterexamples $\CC M(\varphi)$ of an inference formula~$\varphi$
is empty. Hence, we obtain proofs for the \emph{validity} of true inference
rules. On the other hand, if an inference formula is wrong, there must be
a counterexample. A famous theorem in model theory implies that this
counterexample can be chosen algebraic; see~\cite[Proposition~5.2.3]{RealAlgebra}.

\begin{genthm}[Tarski's transfer principle]
The semialgebraic set defined by $\Set{ f_i = 0, g_j \ge 0, h_k \not= 0 }$
is non-empty over some real-closed field if and only if it is non-empty
over every real-closed field, in particular the real closure of~$\BB Q$.
\end{genthm}

The counterexample is a symmetric matrix $\Sigma \in \SymMat_N(\BB K)$
whose entries come from some finite real extension $\BB K$ of $\BB Q$.
By the Primitive element theorem $\BB K = \BB Q(\alpha)$ and all entries
of $\Sigma$ may be represented exactly on a computer as polynomials modulo
the minimal polynomial of $\alpha$. This allows again verification of a
claim about the \emph{invalidity} of an inference formula by off-the-shelf
computer algebra systems. In summary:

\begin{genthm}[Alternatives in Gaussian CI~inference]
If $\varphi$ is a true inference rule for Gaussians, there exists a final
polynomial proof for it with integer coefficients. Otherwise there exists
a counterexample to $\varphi$ with real algebraic coordinates.
\end{genthm}

The Matúš identity is a final polynomial for weak transitivity (and
all other gaussoid axioms).
For four and more Gaussian random variables, there exist inference
rules which are valid but do not follow from the gaussoid axioms.
To be precise, for $n=4$ there are five of them up to symmetry, for example
\cite[Lemma~10, eq.~(20)]{LnenickaMatus}
\begin{align}
  \label{eq:lm20}
  \CI{i,j|k} \wedge \CI{i,k|l} \wedge \CI{i,l|j} \Rightarrow \CI{i,j}.
\end{align}
This inference rule is valid for $4 \times 4$ positive definite matrices.
From its proof one can extract the following fact: the bracket polynomial
\begin{gather}
  \label{eq:lm20f}
  \begin{gathered}
  \apr{ij|\emptyset} \cdot \Big(
    \pr{jk} \apr{jl|\emptyset}^2 \apr{kl|\emptyset}^2 + {} \qquad\quad \\
    \qquad\quad \pr{j} \pr{k}^2 \pr{l} \pr{jl} + \pr{j} \pr{k} \pr{kl}
    \apr{jl|\emptyset}^2
  \Big)
  \end{gathered}
\end{gather}
is in the ideal generated by the assumptions of~\eqref{eq:lm20}. But this
polynomial splits into the desired conclusion $\apr{ij|\emptyset}$ and
another factor which is a sum of non-negative polynomials at least one of
which is positive as a product of principal minors. Hence,
$\apr{ij|\emptyset}$ must vanish whenever the assumptions are satisfied
by a positive definite matrix.

Notably, counterexamples to \eqref{eq:lm20} exist when positive definiteness
is relaxed to allow (indefinite) matrices all of whose principal minors do
not vanish. This condition of \emph{principal regularity} is a natural
substitute for positive definiteness over algebraically closed fields.
An example is this matrix:
\[
  \kbordermatrix{
    &  i &  j &  k &  l \\
  i &  1 &  4 & -2 & -8 \\
  j &  4 &  1 & -2 & -2 \\
  k & -2 & -2 &  1 & \sfrac14 \\
  l & -8 & -2 & \sfrac14 & 1
  }.
\]
This shows that positive definiteness is a crucial feature of our reasoning
task. In particular, we may not find all valid inference rules by working
in an algebraically closed field and basing all computations on ideals only.

\Ueberschriftu{Computation and hardness}

The problem of deciding whether a proposed inference formula $\varphi$ is
valid for all Gaussian distributions over ground set~$N$ reduces the problem
of checking if the model of counterexamples $\CC M(\varphi)$ is empty.
This is a problem in the \emph{existential theory of the reals} (ETR) and
the associated complexity class of all problems which reduce in polynomial
time to ETR is known as $\exists\BB R$; cf.~\cite{ETR}. This is a fundamental
complexity class in computational geometry, polynomial optimization and
statistics. We have thus an upper bound on the complexity of the Gaussian
CI~inference problem.

Unfortunately this upper bound is attained in the complexity-theoretic sense
due to a universality result for Gaussian CI~models. Similar universality
theorems have been known for the realizability of rank-3 matroids
\cite{BokowskiSturmfels}, 4-polytopes \cite{Polytopes} or Nash equilibria
of 3-person games \cite{Datta}. See \cite[Chapter~5]{Dissert} for a more
thorough discussion.

\begin{theorem}[\cite{Dissert}]
The Gaussian CI~inference problem is $\text{co-$\exists\BB R$-complete}$.
\end{theorem}

This theorem is proved by encoding synthetic geometry in the real projective
plane in conditional independence and dependence constraints. Even though
CI~models are at first glance very special semialgebraic sets, they possess
no structure which makes them in general easier to work with than arbitrary
semialgebraic sets.

There exist implementations of quantifier elimination over the real numbers
via cylindrical algebraic decomposition in computer algebra systems such as
Wolfram \TT{Mathematica}. These methods, when they terminate, decide the
emptiness of a semialgebraic set and in the inhabited case they return a
real algebraic number certifying this. Software for finding final polynomials
does not seem to be readily available but there are recent advances in
computing real radical ideals~\cite{Radical}.

An even more fundamental obstacle in the tabulation of all valid inference
rules among five Gaussian random variables is finding reasonable candidate
implications. The following problem should enable an experimental and
data-driven approach:

\begin{problem}
Develop (numerical) software for sampling positive definite points
uniformly from varieties inside $\SymMat_N(\BB R)$.
\end{problem}

Sampling allows to check if a proposed inference formula has any obvious
counterexamples. It can also aid in testing candidates for final polynomials,
such as the left-hand side in \eqref{eq:pappus} or the expression
\eqref{eq:lm20f}, by evaluating these candidates on sufficiently many
samples and checking whether they vanish. Since CI~equations are continuous,
small numerical errors can be tolerated.



\end{multicols}

\end{otherlanguage}

\end{document}

%% file: pappus.tikz

\begin{scope}
  \definecolor{c00a900}{RGB}{0,220,0}
  \definecolor{c0000a9}{RGB}{0,0,220}
  \definecolor{ca90000}{RGB}{220,0,0}

  \path[cm={{1.0,0.0,0.0,1.0,(-66.0,-5.0)}},draw=black,line join=miter,line cap=round,miter limit=3.25,line width=2.543pt] (66.0000,186.0273) -- (642.0000,70.8281);
  \path[cm={{1.0,0.0,0.0,1.0,(-66.0,-5.0)}},draw=black,line join=miter,line cap=round,miter limit=3.25,line width=2.543pt] (66.0000,383.5156) -- (642.0000,325.9141);
  \path[cm={{1.0,0.0,0.0,1.0,(-66.0,-5.0)}},draw=c00a900,line join=miter,line cap=round,miter limit=3.25,line width=2.543pt] (359.4844,457.5703) -- (510.3438,5.0000);
  \path[cm={{1.0,0.0,0.0,1.0,(-66.0,-5.0)}},draw=c00a900,line join=miter,line cap=round,miter limit=3.25,line width=2.543pt] (96.3828,457.5703) -- (583.7656,5.0000);
  \path[cm={{1.0,0.0,0.0,1.0,(-66.0,-5.0)}},draw=c0000a9,line join=miter,line cap=round,miter limit=3.25,line width=2.543pt] (586.1094,457.5703) -- (329.9336,5.0000);
  \path[cm={{1.0,0.0,0.0,1.0,(-66.0,-5.0)}},draw=c0000a9,line join=miter,line cap=round,miter limit=3.25,line width=2.543pt] (118.6094,457.5703) -- (489.5703,5.0000);
  \path[cm={{1.0,0.0,0.0,1.0,(-66.0,-5.0)}},draw=ca90000,line join=miter,line cap=round,miter limit=3.25,line width=2.543pt] (66.0000,47.3203) -- (642.0000,417.6055);
  \path[cm={{1.0,0.0,0.0,1.0,(-66.0,-5.0)}},draw=ca90000,line join=miter,line cap=round,miter limit=3.25,line width=2.543pt] (484.2852,457.5703) -- (107.1445,5.0000);
  \path[cm={{1.0,0.0,0.0,1.0,(-66.0,-5.0)}},draw=black,dash pattern=on 9.22pt off 13.82pt,line join=miter,line cap=round,miter limit=3.25,line width=2.543pt] (66.0000,286.4805) -- (642.0000,163.0508);

  \path[fill=black,even odd rule] (417.0469,99.0000) .. controls (417.0469,97.3945) and (416.4102,95.8594) .. (415.2773,94.7227) .. controls (414.1406,93.5898) and (412.6055,92.9531) .. (411.0000,92.9531) .. controls (409.3945,92.9531) and (407.8594,93.5898) .. (406.7227,94.7227) .. controls (405.5898,95.8594) and (404.9531,97.3945) .. (404.9531,99.0000) .. controls (404.9531,100.6055) and (405.5898,102.1406) .. (406.7227,103.2773) .. controls (407.8594,104.4102) and (409.3945,105.0469) .. (411.0000,105.0469) .. controls (412.6055,105.0469) and (414.1406,104.4102) .. (415.2773,103.2773) .. controls (416.4102,102.1406) and (417.0469,100.6055) .. (417.0469,99.0000) -- cycle(417.0469,99.0000);
  \path[fill=black,even odd rule] (335.0469,115.0000) .. controls (335.0469,113.3945) and (334.4102,111.8594) .. (333.2773,110.7227) .. controls (332.1406,109.5898) and (330.6055,108.9531) .. (329.0000,108.9531) .. controls (327.3945,108.9531) and (325.8594,109.5898) .. (324.7227,110.7227) .. controls (323.5898,111.8594) and (322.9531,113.3945) .. (322.9531,115.0000) .. controls (322.9531,116.6055) and (323.5898,118.1406) .. (324.7227,119.2773) .. controls (325.8594,120.4102) and (327.3945,121.0469) .. (329.0000,121.0469) .. controls (330.6055,121.0469) and (332.1406,120.4102) .. (333.2773,119.2773) .. controls (334.4102,118.1406) and (335.0469,116.6055) .. (335.0469,115.0000) -- cycle(335.0469,115.0000);
  \path[fill=black,even odd rule] (171.0469,148.0000) .. controls (171.0469,146.3945) and (170.4102,144.8594) .. (169.2773,143.7227) .. controls (168.1406,142.5898) and (166.6055,141.9531) .. (165.0000,141.9531) .. controls (163.3945,141.9531) and (161.8594,142.5898) .. (160.7227,143.7227) .. controls (159.5898,144.8594) and (158.9531,146.3945) .. (158.9531,148.0000) .. controls (158.9531,149.6055) and (159.5898,151.1406) .. (160.7227,152.2773) .. controls (161.8594,153.4102) and (163.3945,154.0469) .. (165.0000,154.0469) .. controls (166.6055,154.0469) and (168.1406,153.4102) .. (169.2773,152.2773) .. controls (170.4102,151.1406) and (171.0469,149.6055) .. (171.0469,148.0000) -- cycle(171.0469,148.0000);
  \path[fill=black,even odd rule] (459.0469,333.0000) .. controls (459.0469,331.3945) and (458.4102,329.8594) .. (457.2773,328.7227) .. controls (456.1406,327.5898) and (454.6055,326.9531) .. (453.0000,326.9531) .. controls (451.3945,326.9531) and (449.8594,327.5898) .. (448.7227,328.7227) .. controls (447.5898,329.8594) and (446.9531,331.3945) .. (446.9531,333.0000) .. controls (446.9531,334.6055) and (447.5898,336.1406) .. (448.7227,337.2773) .. controls (449.8594,338.4102) and (451.3945,339.0469) .. (453.0000,339.0469) .. controls (454.6055,339.0469) and (456.1406,338.4102) .. (457.2773,337.2773) .. controls (458.4102,336.1406) and (459.0469,334.6055) .. (459.0469,333.0000) -- cycle(459.0469,333.0000);
  \path[fill=black,even odd rule] (335.0469,346.0000) .. controls (335.0469,344.3945) and (334.4102,342.8594) .. (333.2773,341.7227) .. controls (332.1406,340.5898) and (330.6055,339.9531) .. (329.0000,339.9531) .. controls (327.3945,339.9531) and (325.8594,340.5898) .. (324.7227,341.7227) .. controls (323.5898,342.8594) and (322.9531,344.3945) .. (322.9531,346.0000) .. controls (322.9531,347.6055) and (323.5898,349.1406) .. (324.7227,350.2773) .. controls (325.8594,351.4102) and (327.3945,352.0469) .. (329.0000,352.0469) .. controls (330.6055,352.0469) and (332.1406,351.4102) .. (333.2773,350.2773) .. controls (334.4102,349.1406) and (335.0469,347.6055) .. (335.0469,346.0000) -- cycle(335.0469,346.0000);
  \path[fill=black,even odd rule] (129.0469,366.0000) .. controls (129.0469,364.3945) and (128.4102,362.8594) .. (127.2773,361.7227) .. controls (126.1406,360.5898) and (124.6055,359.9531) .. (123.0000,359.9531) .. controls (121.3945,359.9531) and (119.8594,360.5898) .. (118.7227,361.7227) .. controls (117.5898,362.8594) and (116.9531,364.3945) .. (116.9531,366.0000) .. controls (116.9531,367.6055) and (117.5898,369.1406) .. (118.7227,370.2773) .. controls (119.8594,371.4102) and (121.3945,372.0469) .. (123.0000,372.0469) .. controls (124.6055,372.0469) and (126.1406,371.4102) .. (127.2773,370.2773) .. controls (128.4102,369.1406) and (129.0469,367.6055) .. (129.0469,366.0000) -- cycle(129.0469,366.0000);
  \path[fill=black,even odd rule] (383.0469,201.0000) .. controls (383.0469,199.3945) and (382.4102,197.8594) .. (381.2773,196.7227) .. controls (380.1406,195.5898) and (378.6055,194.9531) .. (377.0000,194.9531) .. controls (375.3945,194.9531) and (373.8594,195.5898) .. (372.7227,196.7227) .. controls (371.5898,197.8594) and (370.9531,199.3945) .. (370.9531,201.0000) .. controls (370.9531,202.6055) and (371.5898,204.1406) .. (372.7227,205.2773) .. controls (373.8594,206.4102) and (375.3945,207.0469) .. (377.0000,207.0469) .. controls (378.6055,207.0469) and (380.1406,206.4102) .. (381.2773,205.2773) .. controls (382.4102,204.1406) and (383.0469,202.6055) .. (383.0469,201.0000) -- cycle(383.0469,201.0000);
  \path[fill=black,even odd rule] (285.0469,222.0000) .. controls (285.0469,220.3945) and (284.4102,218.8594) .. (283.2773,217.7227) .. controls (282.1406,216.5898) and (280.6055,215.9531) .. (279.0000,215.9531) .. controls (277.3945,215.9531) and (275.8594,216.5898) .. (274.7227,217.7227) .. controls (273.5898,218.8594) and (272.9531,220.3945) .. (272.9531,222.0000) .. controls (272.9531,223.6055) and (273.5898,225.1406) .. (274.7227,226.2773) .. controls (275.8594,227.4102) and (277.3945,228.0469) .. (279.0000,228.0469) .. controls (280.6055,228.0469) and (282.1406,227.4102) .. (283.2773,226.2773) .. controls (284.4102,225.1406) and (285.0469,223.6055) .. (285.0469,222.0000) -- cycle(285.0469,222.0000);
    \path[fill=black,even odd rule] (240.0469,231.0000) .. controls (240.0469,229.3945) and (239.4102,227.8594) .. (238.2773,226.7227) .. controls (237.1406,225.5898) and (235.6055,224.9531) .. (234.0000,224.9531) .. controls (232.3945,224.9531) and (230.8594,225.5898) .. (229.7227,226.7227) .. controls (228.5898,227.8594) and (227.9531,229.3945) .. (227.9531,231.0000) .. controls (227.9531,232.6055) and (228.5898,234.1406) .. (229.7227,235.2773) .. controls (230.8594,236.4102) and (232.3945,237.0469) .. (234.0000,237.0469) .. controls (235.6055,237.0469) and (237.1406,236.4102) .. (238.2773,235.2773) .. controls (239.4102,234.1406) and (240.0469,232.6055) .. (240.0469,231.0000) -- cycle(240.0469,231.0000);
\end{scope}